\documentclass{gtart_h}

\def\ifplaintex{\expandafter\ifx\csname documentclass\endcsname\relax}

\def\gtp{{\mathsurround=0pt\it $\cal G\mskip-2mu$eometry \&\ 
$\cal T\!\!$opology $\cal P\!$ublications}}  

\def\Addressesr{\bigskip
{\small \parskip 0pt \leftskip 0pt \rightskip 0pt plus 1fil \def\\{\par}
\sl\theaddress\par
\medskip
\rm Email:\stdspace\tt\theemail\hfill\rm Received:\qua\receiveddate \par}}

\def\recd{{\small Received:\qua\receiveddate\ifx\reviseddate\relax
\else\qquad Revised:\qua\reviseddate\fi\par}} 

\def\lognumber#1{\def\thelognumber{#1}}
\def\volumenumber#1{\def\thevolumenumber{#1}}
\def\volumeyear#1{\def\thevolumeyear{#1}}
\def\papernumber#1{\def\thepapernumber{#1}}
\def\pagenumbers#1#2{\def\startpage{#1}\def\finishpage{#2}}
\def\published#1{\def\publishdate{#1}}

\def\received#1{\def\receiveddate{#1}}

\def\accepted#1{\def\accepteddate{#1}}
\def\asciititle#1{\def\theasciititle{#1}}
\def\covertitle#1{\def\thecovertitle{#1}}

\long\def\asciiabstract#1{\long\def\theasciiabstract{#1}}

\let\\\par\let\thelognumber\relax\let\thevolumenumber\relax
\let\thepapernumber\relax\let\thevolumeyear\relax\let\startpage\relax
\let\finishpage\relax\let\publishdate\relax\let\receiveddate\relax
\let\reviseddate\relax\let\accepteddate\relax\let\theasciititle\relax
\let\thecovertitle\relax\let\theasciiauthors\relax
\let\theasciiabstract\relax

\let\theasciiemail\relax

\ifplaintex
\font\logobig=cmssbx10 scaled 3836
\font\logomed=cmssbx10 scaled 2557
\else
\font\logobig=cmssbx10 scaled 4200
\font\logomed=cmssbx10 scaled 2800
\fi

\long\def\makeagttitle{   
\count0=\startpage
\agt\hfill      
\hbox to 45truept{\vbox to 0pt{\vglue -13truept{\logomed A\kern -.37em{\logobig 
T}\kern -.38em G}\vss}\hss}
\break
{\small Volume \thevolumenumber\ (\thevolumeyear)
\startpage--\finishpage\nl
Published: \publishdate}

\vglue .25truein

{\parskip=0pt\leftskip 0pt plus
1fil\def\\{\par\smallskip}{\Large\bf\thetitle}\par\medskip} \vglue
0.05truein

{\parskip=0pt\leftskip 0pt plus 1fil\def\\{\par}{\sc\theauthors}
\par\medskip}%
 
\vglue 0.03truein 

{\small\leftskip 25truept\rightskip 25truept{\bf Abstract}\stdspace\theabstract

{\bf AMS Classification}\stdspace\theprimaryclass
\ifx\thesecondaryclass\relax\else; \thesecondaryclass\fi\par
{\bf Keywords}\stdspace \thekeywords\par}\vglue 7truept

}   

\ifplaintex
\font\phead=cmsl9 scaled 950
\font\pnum=cmbx10 scaled 913
\font\pfoot=cmsl9 scaled 950
\headline{\vbox to 0pt{\vskip -4.5mm\line{\small\phead\ifnum
\count0=\startpage ISSN 1472-2739 (on-line) 1472-2747 (printed)
\hfill {\pnum\folio}\else\ifodd\count0\def\\{ }%
\ifx\theshorttitle\relax\thetitle\else\theshorttitle\fi\hfill{\pnum\folio}
\else\def\\{ and }{\pnum\folio}\hfill\ifx\theshortauthors\relax\theauthors
\else\theshortauthors\fi\fi\fi}\vss}}
\footline{\vbox to 0pt{\vglue 0mm\line{\small\pfoot\ifnum\count0=\startpage
\copyright\ \gtp\hfill\else
\agt, Volume \thevolumenumber\ (\thevolumeyear)\hfill\fi}\vss}}
\else
\font=cmsl9 scaled 1050
\font\lnum=cmbx10 
\font=cmsl9 scaled 1050
\makeatletter
\def\@oddhead{{\small\lhead\ifnum\count0=\startpage ISSN 1472-2739 
(on-line) 1472-2747 (printed)\hfill {\lnum\number\count0}\else\ifodd\count0
\def\\{ }\ifx\theshorttitle\relax \thetitle \else\theshorttitle\fi\hfill
{\lnum\number\count0}\else\def\\{ and }{\lnum\number\count0}
\hfill\ifx\theshortauthors\relax 
\theauthors\else\theshortauthors\fi\fi\fi}}\def\@evenhead{\@oddhead}
\def\@oddfoot{\small\lfoot\ifnum\count0=\startpage\copyright\ \gtp\hfill\else
\agt, Volume \thevolumenumber\ (\thevolumeyear)\hfill\fi}
\def\@evenfoot{\@oddfoot}
\makeatother
\fi
\let\maketitlepage\makeagttitle
\let\makeshorttitle\maketitlepage
\let\maketitle\maketitlepage

\newwrite\gtoutfile
\long\gdef\makeheadfile{  
{\def\\{, }\def\s{ }
\immediate\openout\gtoutfile head.xxx
\immediate\write\gtoutfile{Proxy-for: \ifx\theasciiauthors\relax
\theauthors\else\theasciiauthors\fi\s<\ifx\theasciiemail\relax\theemail\else\theasciiemail\fi>}
\immediate\write\gtoutfile{\noexpand\\}
\immediate\write\gtoutfile{Authors: \ifx\theasciiauthors\relax
\theauthors\else\theasciiauthors\fi}
{\def\\{ }\immediate\write\gtoutfile{Title: \ifx\theasciititle\relax
\thetitle\else\theasciititle\fi}}
\immediate\write\gtoutfile{Subj-class: GT or SG, GR etc}
\immediate\write\gtoutfile{MSC-class: \theprimaryclass\ifx\thesecondaryclass\relax\else, \thesecondaryclass\fi}
\immediate\write\gtoutfile{Journal-ref: Algebr. Geom. Topol. \thevolumenumber\s
(\thevolumeyear) \startpage-\finishpage}
\immediate\write\gtoutfile{Comments: Published by Algebraic and
Geometric Topology at}
\immediate\write\gtoutfile{\s\s\s  http://www.maths.warwick.ac.uk/agt/AGTVol\thevolumenumber/agt-\thevolumenumber-\thepapernumber.abs.html}
\immediate\write\gtoutfile{\noexpand\\}
\immediate\write\gtoutfile{}
\ifx\theasciiabstract\relax
\immediate\write\gtoutfile{\theabstract}\else
\immediate\write\gtoutfile{\theasciiabstract}\fi
\immediate\write\gtoutfile{}
\immediate\write\gtoutfile{\noexpand\\}
\immediate\write\gtoutfile{}
\immediate\closeout\gtoutfile}}  

\def\maketitlepage{\makeagttitle\makeheadfile}
\let\makeshorttitle\maketitlepage
\let\maketitle\maketitlepage

\lognumber{10}
\volumenumber{5}
\volumeyear{2005}
\papernumber{10}
\pagenumbers{183}{205}
\received{31 May 2004} 
\accepted{21 September 2004}
\published{23 March 2005}

\usepackage{amssymb,amsmath,amscd}

\newtheorem{thm}{Theorem}[section]    
\newtheorem{lem}[thm]{Lemma}          
\newtheorem{cor}[thm]{Corollary}
\newtheorem{prop}[thm]{Proposition}

\newtheorem{condition}[thm]{Condition}

\theoremstyle{definition}
  
\newtheorem{rem}[thm]{Remark}             
\newtheorem{example}[thm]{Example}

\newcommand{\nbiga}{\mathcal{A}}
\newcommand{\nbigp}{\mathcal{P}}
\newcommand{\nbigq}{\mathcal{Q}}

\newcommand{\seisuu}{\mathbb{Z}}
\newcommand{\Fin}{\mathbb{F}}
\newcommand{\Fq}{\mathbb{F}_q}
\newcommand{\lrarr}{\longrightarrow}

\def\Image{\mathop{\rm Im}\nolimits}
\def\Ker{\mathop{\rm Ker}\nolimits}
\def\modulo{\mathop{\rm mod}\nolimits}

\newcommand{\gyakuzou}{\delta^{-1}(\Image D_3^{(s)})}

\newcommand{\defn}{\emph}

\begin{document}
\title{The 3--cocycles of the Alexander
quandles\\${\mathbb F}_q[T]/(T{-}\omega)$}

\asciititle{The 3-cocycles of the Alexander quandles F_q[T]/(T-omega)}

\covertitle{The 3--cocycles of the Alexander
quandles\\${\noexpand\bf F}_q[T]/(T{-}\omega)$}
\authors{Takuro Mochizuki}                  
\address{Department of Mathematics, Kyoto University\\Kyoto, 606--8502, Japan}
\email{takuro@math.kyoto-u.ac.jp}
\begin{abstract} 
We determine the third cohomology of Alexander quandles of the
form $\Fq[T]/(T-\omega)$, where $\Fq$ denotes the finite field of order $q$
and $\omega$ is an element of $\Fq$ which is neither $0$ nor $1$.
As a result, we obtain many concrete examples of non-trivial $3$--cocycles.
\end{abstract}

\asciiabstract{%
We determine the third cohomology of Alexander quandles of the form
F_q[T]/(T-omega), where F_q denotes the finite field of order q and
omega is an element of F-q which is neither 0 nor 1.  As a result, we
obtain many concrete examples of non-trivial 3-cocycles.}

\primaryclass{18H40} 
\secondaryclass{55A25, 57Q45}              
\keywords{Quandle, cohomology, knot} 

\maketitle 

\section{Introduction}
\label{section;05.2.19.1}

\subsection{Aims}

Interest in quandles has been growing recently, particularly because of
their applications to the study of classical knots and 2--knots.
A \defn{quandle} is a set $X$ equipped with a binary operation $\ast\co X
\times X \to X$
satisfying the following conditions:
\begin{description}
\item[(Idempotency)]
For any $a\in X$, $a\ast a=a$.
\item[(Right-Invertibility)]
 For any $a,b\in X$, there exists a unique $c\in X$
 such that $a=c\ast b$.
\item[(Self-Distributivity)]
The identity $(a\ast b)\ast c=(a\ast c)\ast (b\ast c)$ holds for all
$a,b,c\in X$.
\end{description}
Quandle cohomology $H^{\ast}(X,A)$ is defined
for any quandle $X$ and any abelian group $A$, and may be used (see
\cite{q1,q2,rs} for details) to construct isotopy invariants of
classical knots and links, and also of higher-dimensional embeddings.  In
particular, such invariants obtained from 3--cocycles play an interesting
r\^ole in the study of 2--knots (see \cite{ss1,ss2} for example).
However, there are not very many concrete examples of nontrivial
3--cocycles, and so it would be useful to find a systematic method of
constructing nontrivial 3--cocycles with calculable forms for some classes
of quandles.

In this paper, we discuss the third cohomology group of Alexander
quandles.
Let $R$ be a commutative ring with a unit element,
and let $M$ be an $R$--module.
For any invertible element $\omega$ of $R$,
a binary operation $\ast\co M \times M \to M$ may be defined by
$$a\ast b=\omega\cdot a+(1-\omega)\cdot b.$$
It is easy to check that $(M,\ast)$ satisfies the three quandle axioms,
and we call a quandle of this type an \defn{Alexander quandle}.
In this paper, we restrict ourselves to the case $R=M=\Fq$, where
$q$ is a power of a prime $p$ and $\Fq$ denotes a finite field
of order $q$, and denote the resulting quandle
$\Fq[T]/(T-\omega)$.  We ignore the case $\omega=1$, which
yields a trivial quandle, and the case $\omega=0$, which is forbidden by
the right-invertibility axiom.

Quandle cohomology groups are well understood in the case where $A$ is a
field of characteristic $0$ (see \cite{eg,m}), and do not typically give
rise to interesting cocycles.  However, we may expect interesting examples
to arise from the case where $A$ is a field of positive characteristic.
For example, the third quandle cohomology group
$H^3\bigl(\Fin_p[T]/(T-\omega),\Fin_p\bigr)$
was calculated in \cite{m}, and the case $\omega=-1$ gives rise to a
particular nontrivial 3--cocycle which has been used by Satoh and Shima in
their study of 2--knots.

By generalizing our previous work, we will determine the third quandle
cohomology group $H^3\bigl(\Fin_q[T]/(T-\omega),A\bigr)$ where $A$
is an algebraic closure $k$ of $\Fq$ (this includes the case where
$A$ is a field of characteristic $p$).  In doing so (see subsections
\ref{subsection;05.2.19.55} and \ref{subsection;9.25.2}), we obtain many
examples of nontrivial 3--cocycles, which we hope will be useful in the
study of 2--knots.

\subsection{Outline}
In subsection \ref{subsection;05.2.19.50},
we recall the definition of the quandle cohomology groups
and explain our description of the cocycles, which is slightly
different to the description in \cite{m}.
In subsection \ref{subsection;05.2.19.51}, we give some concrete examples.
Subsection \ref{subsection;05.2.19.55} contains Theorem~\ref{thm;9.16.1},
the main result of this paper, which we apply to certain quandles in
subsection \ref{subsection;9.25.2}.
The proof of the theorem is given in Section~\ref{section;9.25.3}.

\subsection{Acknowledgements}
The author is grateful to the referee for his comments, and for drawing
attention to a paper of Rourke and Sanderson \cite{rs}.
The author thanks Mikiya Masuda, Yoshifumi Tsuchimoto and Akira Ishii
for their encouragement,
Nicholas Jackson for polishing the English translation,
the Japan Society for the Promotion of Science,
and the Sumitomo Foundation, for financial support, and the Institute of
Advanced Study (where this paper was written) for their kind hospitality.
This material is based upon work supported by National Science Foundation
grant DMS--9729992, although any opinions, findings, conclusions or
recommendations expressed herein are the responsibility of the author.

\section{The main result} \label{section;9.25.1}

\subsection{Preliminaries}
\label{subsection;05.2.19.50}

\subsubsection{Quandle cohomology groups}
\label{subsection;05.2.19.20}

Let $(X,\ast)$ be a quandle, and let $A$ be an abelian group.
We define a complex $C^\ast(X,A)$ with cochain groups
$$C^{n}\bigl(X,A\bigr):=
  \bigl\{f\co X^n\to A\,\big|\,
  f(x_1,\ldots,x_n)=0\,\,\mbox{\rm when }
  x_i=x_{i+1}\,\,\mbox{ for some }i \bigr\}$$
and differential $\delta\co C^n(X,A)\lrarr C^{n+1}(X,A)$ defined
as follows:
\begin{multline*}
 \delta(f)(x_1,\ldots,x_{n+1}):=
 \sum_{i=1}^{n+1}
 (-1)^{i-1} f\bigl(
 x_1\ast x_i,\ldots, x_{i-1}\ast x_i,
 x_{i+1},\ldots,x_{n+1}
 \bigr) \\
-\sum_{i=1}^{n+1}
 (-1)^{i-1} f\bigl(x_1,\ldots,x_{i-1},x_{i+1},\ldots,x_{n+1}\bigr)
\end{multline*}
The cohomology of this complex is denoted $H^*(X,A)$, and called the
\defn{quandle cohomology} of $X$ (with coefficient group $A$).

Let $q$ be a power of a prime $p$, let $\Fq$ denote a finite field
of order $q$,
and let $\omega$ be a non-zero element of $\Fq$.
We wish to calculate
$H^{3}\bigl(\Fq[T]/(T-\omega),k\bigr)$, where $\Fq[T]/(T-\omega)$ is the
Alexander quandle discussed in section \ref{section;05.2.19.1}, and $k$ is
an algebraic closure of $\Fq$.  It is obvious that the differential map is
trivial in the case $\omega=1$, so we will consider only the cases $\omega\neq
0,1$.

\subsubsection{The quandle complex}
\label{subsubsection;05.2.19.2}

Let $k$ be a field, and $\omega$ an element of $k$ which is neither 0 nor
1.  (Later $k$ will denote an algebraic closure of $\Fq$, but for the
moment we consider arbitrary fields.)
Let $k[U_1,\ldots,U_n]$ be the polynomial ring
over $k$ with $n$ variables $U_1,\ldots,U_{n}$, and
set $\Omega_{n-1}:=\prod_{i=1}^{n-1}U_i$
and $C^n:=\Omega_{n-1}\cdot k[U_1,\ldots,U_n]$.

For any element $f\in C^n$,
we define $\delta(f)\in C^{n+1}$ as follows:
\begin{multline*}
 \delta(f)\bigl(U_1,\ldots,U_{n+1}\bigr):= \\
 \sum_{i=1}^{n}
 (-1)^{i-1}
 f\bigl(\omega\cdot  U_1,\ldots,
 \omega\cdot  U_{i-1},\omega\cdot  U_{i}+U_{i+1},
 U_{i+2},\ldots,U_{n+1}\bigr) \\
-\sum_{i=1}^{n-1}
 (-1)^{i-1}
 f\bigl(U_1,\ldots,U_{i-1},U_{i}+U_{i+1},U_{i+2},
  \ldots,U_{n+1}\bigr)
\end{multline*}
We thus obtain a homomorphism $\delta\co C^n\to C^{n+1}$.  A routine
calculation verifies that $\delta\circ\delta=0$, and so we have a
complex $C^{\ast}=\big(\bigoplus_{n=1}^{\infty}C^n,\delta\big)$, which we
call the \defn{quandle complex associated with $k$ and $\omega$}.

\begin{rem}
The complex $C^\ast$ was discussed in \cite{m} in the case where $k$ is a
field of characteristic $0$, and was shown to be acyclic.  The above
definition looks slightly different, due to a different choice of
coordinates in $k^n$.
\end{rem}

\subsubsection{A convenient description
of the complex $C^{\ast}\bigl(\Fin_q[T]/(T-\omega),k\bigr)$}
\label{subsection;9.25.5}

Let $\Fin_q$ denote a finite field of order $q$.
In the following,
$k$ is an algebraic closure of $\Fin_q$,
and $\omega$ an element of $\Fin_q$ such that $\omega\neq 0,1$.
Let $C^{\ast}$ be the complex described in
subsection \ref{subsubsection;05.2.19.2}.  Then an element
$f=f(U_1,\ldots,U_n)$ of $C^n=k[U_1,\ldots,U_n]$
induces a $k$--valued function on $\Fin_q^n$, given by
\[
 \Fin_q^n
 \ni
 \bigl(x_1,\ldots,x_{n}\bigr)\mapsto
 f\bigl(x_1-x_2,x_2-x_3,\ldots,x_{n-1}-x_n,x_n\bigr)
 \in k,
\]
and so we obtain a map
$\varphi\co C^n\to C^n\bigl(\Fin_q[T]/(T-\omega),k\bigr)$ for any $n$.
The following lemma can be checked by a direct calculation.
\begin{lem}
The homomorphism $\varphi$ is compatible with the differentials (that is,
$\varphi\circ\delta=\delta\circ\varphi$) and so
we obtain a morphism of the cochain complexes
$C^{\ast}\to C^{\ast}\bigl(\Fin_q[T]/(T-\omega),k\bigr)$.
\end{lem}
We now define
\[
 C^n(q):=
 \Bigl\{
 \sum a_{i_1,\ldots,i_n}\cdot 
 U_1^{i_1}\cdots U_{n}^{i_{n}}
  \in C^n\,\Big|\,
 0\leq i_j\leq q-1
 \Bigr\}.
\]
It is easy to check that $\delta\bigl(C^n(q)\bigr)\subset C^{n+1}(q)$.
Thus, we obtain the subcomplex
$C^{\ast}(q):=\bigl(\bigoplus C^n(q),\delta\bigr)$.
The following lemma can be checked easily.
\begin{lem}
The induced morphism
$C^{\ast}(q)\to C^{\ast}\bigl(\Fin_q[T]/(T-\omega),k\bigr)$
is an isomorphism.
\end{lem}

\subsection{Some 3--cocycles of the complex $C^{\ast}$}
\label{subsection;05.2.19.51}

We now give some concrete examples of 3--cocycles in the complex
$C^{\ast}$.  The last variable $U_n$ for any element $f(U_1,\ldots,U_n)$
of $C^n$ is denoted by $T_n$ in the following argument, and will be
useful for specific calculations.

\subsubsection{The cocycles $\Psi, E_0$ and $E_1$}
\label{subsection;04.5.29.10}

For positive integers $a$ and $b$ set $\mu_a(x,y):=(x+y)^a-x^a-y^a$ and
define the polynomial $\Psi(a,b)\in k[U_1,U_2,T_3]$ as follows:
\begin{multline*}
 \Psi(a,b):=
 \Bigl(
 \mu_a\bigl(\omega\cdot  U_1,U_2\bigr)
-\mu_a\bigl(U_1,U_2\bigr)
 \Bigr)\cdot  T_3^b \\
 =
 \Bigl(
 \bigl(\omega\cdot  U_1+U_2\bigr)^a
  -\bigl(U_1+U_2\bigr)^a+\bigl(1-\omega^a\bigr)\cdot  U_1^a
 \Bigr)\cdot  T_3^b
\end{multline*}

\begin{lem}\label{lem;04.5.26.1}
If $\omega^{a+p^s}=1$, then $\Psi(a,p^s)$ is a quandle 3--cocycle.
\end{lem}
\begin{proof}
Set $h(U_1,T_2):=\mu_a(\omega\cdot  U_1,T_2)-\mu_a(U_1,T_2)$.
Then from the relation
$\delta(T_1^a)=(\omega\cdot  U_1+T_2)^a-(U_1+T_2)^a$,
we see that
$h(U_1,T_2)=\delta(T_1^a)+(1-\omega^a)\cdot  U_1^a$, and so
\[
\delta(h)(U_1,U_2,T_3)=
(1-\omega^a)\cdot \delta(U_1^a)=
(1-\omega^a)\cdot 
h(U_1,U_2).
\]
Then
\begin{multline*}
\delta\bigl(\Psi(a,p^s)\bigr)=
(1-\omega^a)\cdot  h(U_1,U_2)\cdot  T_4^{p^s} \\
+\Bigl(\omega^a\cdot  h(U_1,U_2)\cdot  (\omega\cdot  U_3+T_4)^{p^s}
 -h(U_1,U_2)\cdot  (U_3+T_4)^{p^s}\Bigr)=0,
\end{multline*}
so $\Psi(a,p^s)$ is a quandle cocycle.
\end{proof}

We now introduce the following polynomials
with $\seisuu/p\seisuu$ coefficients:
\[
 \chi(x,y):=
  \sum_{i=1}^{p-1}(-1)^{i-1}\cdot  i^{-1}\cdot  x^{p-i}\cdot  y^i\equiv
\frac{1}{p}\Bigl((x+y)^p-x^p-y^p\Bigr)\quad(\modulo p)
\]
For positive integers $a$ and $b$,
we define the polynomial $E_0\bigl(a\cdot  p,b\bigr)(U_1,U_2,T_3)$
to be:
\[
 E_0\bigl(a\cdot  p,b\bigr):=
 \Bigl(
 \chi(\omega\cdot  U_1,U_2)
-\chi(U_1,U_2)
 \Bigr)^a\cdot  T_3^b
\]
\begin{lem}
If $\omega^{p^{s}+p^h}=1$ and $s>0$, then $E_0(p^s,p^h)$ is a quandle 3--cocycle.
\end{lem}
\begin{proof}
This can be verified by an argument similar to that used in the  proof of 
Lemma~\ref{lem;04.5.26.1}.
\end{proof}
For positive integers $a$ and $b$,
we define the polynomial
$E_1\bigl(a,b\cdot  p\bigr)\in k[U_1,U_2,T_3]$
as
\[
 E_1\bigl(a,b \cdot   p\bigr):=
 U_1^a\cdot 
 \Bigl(
 \chi\bigl(U_2,T_3\bigr)-\chi\bigl(U_2,\omega^{-1}\cdot  T_3\bigr)
 \Bigr)^{b}.
\]
\begin{lem}
If $\omega^{p^t+p^s}=1$ and $s>0$,
then $E_1(p^t,p^s)$ is a quandle 3--cocycle.
\end{lem}
\begin{proof}
Set $h(U_2,T_3):=\chi(U_2,T_3)-\chi(U_2,\omega^{-1}\cdot  T_3)$.
Then
\begin{multline*}
\delta\bigl(E_1(p^t,p^s)\bigr)=
\Bigl(
   \bigl(\omega\cdot  U_1+U_2\bigr)-\bigl(U_1+U_2\bigr)
 \Bigr)^{p^t}\cdot  h(U_3,T_4)^{p^{s-1}} \\
-U_1^{p^t}\cdot \Bigl(
  \omega^{p^t}h\bigl(\omega\cdot  U_2+U_3,T_4\bigr)^{p^{s-1}}
 -h\bigl(U_2+U_3,T_4\bigr)^{p^{s-1}}
 \Bigr)\\
+U_1^{p^t}\cdot \Bigl(
 \omega^{p^t}
  h\bigl(\omega\cdot  U_2,\,\,\omega\cdot  U_3+T_4\bigr)^{p^{s-1}}
-h\bigl(U_2,U_3+T_4\bigr)^{p^{s-1}}
 \Bigr).
\end{multline*}
By using the relation $\omega^{p^t}=\omega^{-p^s}$,
the right hand side can be rewritten
\begin{multline*}
 U_1^{p^t}
\cdot 
\Bigl(
 -\bigl(1-\omega^{-p}\bigr)\cdot  h\bigl(U_3,T_4\bigr)
-\omega^{-p}\cdot  h\bigl(\omega U_2+U_3,T_4\bigr)
+h\bigl(U_2+U_3,T_4\bigr)\\
+\omega^{-p}\cdot  h\bigl(\omega U_2,\omega U_3+T_4\bigr)
-h\bigl(U_2,U_3+T_4\bigr)
\Bigr)^{p^{s-1}}.
\end{multline*}
It can be directly shown that this expression is zero, and so
$E_1(p^t,p^s)$ is a quandle 3--cocycle.
\end{proof}

\subsubsection{The set $\nbigq$ and the cocycles $F$ and $\Gamma$}
\label{subsection;04.5.29.11}

In the following, let $q_i$ be powers of the prime $p$.
For any non-negative integers $a,b,c$ and $d$, we define  polynomials
$F(a,b,c)\in k[U_1,U_2,T_3]$ and $G(a,b,c,d)\in k[U_1,U_2,U_3,T_4]$ as
\begin{align*}
 F(a,b,c)&:= U_1^a\cdot  U_2^b\cdot  T_3^c,\\
 G(a,b,c,d)&:=U_1^a\cdot  U_2^b\cdot  U_3^c\cdot  T_4^d.
\end{align*}
The following lemma is helpful for our later calculations.
\begin{lem} \label{lem;10.20.1}
We have the following identities.
\[
 \begin{array}{rl}
 \delta\bigl(F(q_1,q_2,q_3)\bigr)=&
 \bigl(\omega^{q_1+q_2+q_3}-1\bigr)\cdot  G(q_1,q_2,q_3,0)\\
\\
 \delta\bigl(F(q_1+q_2,q_3,q_4)\bigr)=&
 ~~\bigl(\omega^{q_1}-1\bigr)\cdot  G(q_1,q_2,q_3,q_4) \\
 &
+\bigl(\omega^{q_2}-1\bigr)\cdot  G(q_2,q_1,q_3,q_4)\\
 &
+\bigl(\omega^{q_1+q_2+q_3+q_4}-1\bigr)\cdot  G(q_1+q_2,q_3,q_4,0)\\
\\
\delta\bigl(F(q_1,q_2+q_3,q_4)\bigr)=&
-\bigl(\omega^{q_1+q_2}-1\bigr)\cdot  G(q_1,q_2,q_3,q_4)\\
 &
-\bigl(\omega^{q_1+q_3}-1\bigr)\cdot  G(q_1,q_3,q_2,q_4)\\
 &
+\bigl(\omega^{q_1+q_2+q_3+q_4}-1\bigr)\cdot  G(q_1,q_2+q_3,q_4,0)\\
\\
 \delta\bigl(F(q_1,q_2,q_3+q_4)\bigr)=&
 ~~\bigl(\omega^{q_1+q_2+q_3}-1\bigr)\cdot  G(q_1,q_2,q_3,q_4) \\
&
+\bigl(\omega^{q_1+q_2+q_4}-1\bigr)\cdot  G(q_1,q_2,q_4,q_3)\\
 &
+\bigl(\omega^{q_1+q_2+q_3+q_4}-1\bigr)\cdot  G(q_1,q_2,q_3+q_4)
 \end{array}
\]
\end{lem}
\begin{proof}
These identities may be verified by direct calculation.
\end{proof}

\begin{cor}
Let $q_1,q_2$ and $q_3$ be powers of a prime $p$.  Then
\begin{enumerate}
\item $F(q_1,q_2,q_3)$ is a quandle 3--cocycle if $\omega^{q_1+q_2+q_3}=1$,
and
\item $F(q_1,q_2,0)$ is a quandle 3--cocycle if $\omega^{q_1+q_2}=1$.
\end{enumerate}
\end{cor}

\label{subsection;9.25.11}

Let $\nbigq$ denote the set of 
quadruples $(q_1,q_2,q_3,q_4)$
satisfying the following conditions:
\begin{condition} \label{condition;04.5.28.100} \mbox{{}}
{\rm \begin{itemize}
\item
$q_2\leq q_3$, $q_1<q_3$, $q_2<q_4$, and
$\omega^{q_1+q_3}=\omega^{q_2+q_4}=1$.
\item One of the following holds:
\begin{description}
\item[Case 1] $\omega^{q_1+q_2}=1$.
\item[Case 2] $\omega^{q_1+q_2}\neq 1$ and $q_3>q_4$.
\item[Case 3] ($p\neq 2$) $\omega^{q_1+q_2}\neq 1$ and $q_3=q_4$.
\item[Case 4] ($p\neq 2$) $\omega^{q_1+q_2}\neq 1$, $q_2\leq q_1<q_3<q_4$,
and $\omega^{q_1}=\omega^{q_2}$.
\item[Case 5] ($p=2$) $\omega^{q_1+q_2}\neq 1$, $q_2< q_1<q_3<q_4$,
and $\omega^{q_1}=\omega^{q_2}$.
\end{description}
\end{itemize}}
\end{condition}

The polynomial $\Gamma(q_1,q_2,q_3,q_4)$ is defined 
for any element $(q_1,q_2,q_3,q_4)$ of $\nbigq$ as follows:\\
{\bf Case 1}
$$\Gamma(q_1,q_2,q_3,q_4):= F(q_1,q_2+q_3,q_4).$$
{\bf Case 2}
\begin{multline*}
 \Gamma(q_1,q_2,q_3,q_4):=
 F(q_1,q_2+q_3,q_4)
-F(q_2,q_1+q_4,q_3) \\
-\bigl(\omega^{q_2}-1\bigr)^{-1}\cdot 
 \bigl(1-\omega^{q_1+q_2}\bigr)
\cdot 
 \Bigl(
 F(q_1,q_2,q_3+q_4)-F(q_1+q_2,q_4,q_3)
 \Bigr).
\end{multline*}
{\bf Case 3}
$$\Gamma(q_1,q_2,q_3,q_4):=
 F(q_1,q_2+q_3,q_4)-2^{-1}\cdot  (1-\omega^{-q_3})\cdot 
F(q_1,q_2,q_3+q_4).$$

{\bf Case 4 and Case 5}
\begin{multline*}
 \Gamma(q_1,q_2,q_3,q_4):=
 F(q_1,q_2+q_3,q_4)+F(q_2,q_1+q_3,q_4) \\
-(\omega^{q_1}-1)^{-1}\cdot(1-\omega^{2q_1})\cdot
 F(q_1+q_2,q_3,q_4).
\end{multline*}
The next lemma follows from Lemma~\ref{lem;10.20.1} together with a direct
calculation.
\begin{lem}
The polynomials $\Gamma(q_1,q_2,q_3,q_4)$ are quandle 3--cocycles
for any quadruple $(q_1,q_2,q_3,q_4)\in \nbigq$.
\end{lem}

We define
$\nbigq(q):= \left\{(q_1,q_2,q_3,q_4)\in\nbigq ~|~ q_i<q\right\}$
and, for $d$ a positive integer,
$\nbigq_d(q):= \left\{(q_1,q_2,q_3,q_4)\in\nbigq(q) ~|~ \sum q_i=d\right\}.$
\subsection{Statement of the main theorem}
\label{subsection;05.2.19.55}
As before let $q_i$ denote a power of a prime $p$, and define:
\begin{equation} \label{eq;04.5.26.10}
\begin{aligned} 
 I(q):=&
 \left\{F(q_1,q_2,q_3) ~|~ \omega^{q_1+q_2+q_3}=1,~q_1<q_2<q_3<q\right\} \\
 &\cup
 \left\{
 F(q_1,q_2,0) ~|~ \omega^{q_1+q_2}=1,~q_1<q_2<q
 \right\}\\
 &
\cup
 \left\{\Psi(a,q_1)~\left|~
 \begin{array}{l}
 \omega^{a+q_1}=1,~
 0<a<q,~1<q_1<q,\\
 a\not\equiv 0~(\modulo q_1),~
 \mbox{$a$ is not a power of $p$}
 \end{array}
 \right.\right\}
 \\
 &
\cup
 \left\{
 E_0(p\cdot q_1,q_2) ~|~
 \omega^{p\cdot q_1+q_2}=1,~
 q_1<q_2<q
 \right\} 
 \\
 &\cup
 \left\{
 E_1(q_1,p\cdot q_2) ~|~
\omega^{q_1+p\cdot q_2}=1,~
 q_1\leq q_2<q
 \right\} \\
 &
\cup
 \left\{
 \Gamma(q_1,q_2,q_3,q_4) ~|~
 (q_1,q_2,q_3,q_4)\in\nbigq(q)
 \right\}
\end{aligned}
\end{equation}
Let $H^3(q)$ denote the subspace of $C^3(q)$
generated by $I(q)$.
The following theorem is the main result
of this paper, and will be proved in Section \ref{section;9.25.3}.
\begin{thm} \label{thm;9.16.1}
The natural map
$H^3(q)\to H^3\bigl(C^{\ast}(q)\bigr)$ is an isomorphism.
\end{thm}

\begin{rem}
It will also turn out that
the cocycles given in {\rm(\ref{eq;04.5.26.10})}
are linearly independent, and hence
form a basis for $H^3(q)$.
\end{rem}

\subsection{Examples} \label{subsection;9.25.2}

\subsubsection{The case $\omega=-1$}
Let $p$ be an odd prime, and let $\omega=-1$.  Then:
\begin{itemize}
\item $\omega^{q_1+q_2+q_3}\neq 1$, so we have no
quandle $3$--cocycles of the form $F(q_1,q_2,q_3)$.
\item $\omega^{a+q_1}=-\omega^a$, and so
the identity $\omega^{a+q_1}=1$ implies that $a$ is odd.
\item $\omega^{q_1+q_2}=1$ for any powers $q_i$ of $p$.
Hence the polynomials $F(q_1,q_2,0)$,
$E_0\bigl(p\cdot q_1,q_2\bigr)$ and
$E_1\bigl(q_1,p\cdot q_2\bigr)$ are quandle $3$--cocycles.
In addition,
$\nbigq(q)=\left\{(q_1,q_2,q_3,q_4) ~|~
 q_2\leq q_3,\,q_1<q_3,\,q_2<q_4\right\}$,
and $\omega^{q_1+q_2}=1$ for any
$(q_1,q_2,q_3,q_4)\in \nbigq(q)$.
\end{itemize}
Thus we obtain the following $3$--cocycles,
which form a basis for the cohomology group
$H^3\bigl(\Fq[T]/(T+1),k\bigr)$:
\begin{equation} \label{eq;05.2.19.10}
\begin{aligned}
 & \left\{
 F(q_1,q_2,0)\,\big|\,0<q_1<q_2<q \right\}\\
\cup &
 \left\{
 E_0\bigl(p\cdot q_1,q_2\bigr)~|~q_1<q_2<q
 \right\} \\
\cup &
 \left\{
 E_1\bigl(q_1,p\cdot q_2\bigr)~|~q_1<q_2<q
 \right\}\\
\cup &
 \left\{
 \Psi(a,q_1)~\left|
 \begin{array}{l}
 a~\mbox{odd,}~
 0<a<q,~
 q_1<q,~\\
 a\not\equiv 0~(\modulo q_1),~
 \mbox{$a$ is not a power of $p$}
 \end{array}
 \right.\right\}\\
\cup &
 \left\{
 F(q_1,q_2+q_3,q_4)~|~
 q_2\leq q_3,\,\,q_1<q_3,~q_2<q_4,~
 q_i<q
 \right\}
\end{aligned}
\end{equation}
If $q=p^2$ then the
basis in (\ref{eq;05.2.19.10}) is:
\begin{align*}
 &\left\{
 F(1,p),~ E_0(p,p),~ E_1(1,p),~ E_1(1,p^2),~ E_1(p,p^2),~ F(1,p+1,p) \right\} \\
\cup~&
 \left\{
 \Psi(a,p)~|~ a~\mbox{odd},~ a<p^2,~ a\not\equiv 0~(\modulo p),~ a\neq 1
 \right\}
\end{align*}
If $q=p$ then the basis in (\ref{eq;05.2.19.10}) is simply $\{E_1(1,p)\}$,
and so $H^3\bigl(\Fin_p[T]/(T+1),k\bigr)$ is $1$--dimensional,
as previously noted in \cite{m}.

\subsubsection{Some other examples}
\begin{example} 
If $\Fq=\seisuu_2[\omega]/(1+\omega+\omega^2)$, then $q=4=2^2$,
and the order of $\omega$ is $3=2+1$.  Then:
\begin{itemize}
\item We have no triples $(q_1,q_2,q_3)$
  of powers of $2$ satisfying $q_1<q_2<q_3<2^2$.
\item If a pair $(q_1,q_2)$ of powers of $2$ satisfies $q_1<q_2<2^2$,
then we have $q_1=1$ and $q_2=2$.
In this case, $\omega^{q_1+q_2}=\omega^3=1$, and so we have a cocycle
$F(1,2,0)$.
\item The identity $\omega^{a+2}=1$ implies $a\equiv 1~(\modulo 3)$, and
$a$ cannot be a power of $2$ because of the definition of $I(q)$ in
\eqref{eq;04.5.26.10}, hence cocycles of the form $\Psi(b\cdot p^t,p^s)$
do not occur.
\item If $(q_1,q_2)=(1,2)$,
then $\omega^{2\cdot q_1+q_2}=\omega^4\neq 1$.
\item If $(q_1,q_2)=(1,2)$,
then $\omega^{q_1+2\cdot q_2}=\omega^5\neq 1$.
On the other hand,
if $(q_1,q_2)=(1,1)$ or $(2,2)$,
then $\omega^{q_1+2\cdot q_2}=\omega^{3q_1}=1$, and so
we have the cocycles $E_1(1,2)$ and $E_1(2,4)$.
\item The set $\nbigq(q)$ is empty.
\end{itemize}
Thus we have cocycles
$$\bigl\{ F(1,2),\,\,E_1(1,2),\,\,E_1(2,4) \bigr\}$$
which form a basis for the cohomology group $H^3\bigl(\Fq[T]/(T+1),k\bigr)$.
\end{example}

\begin{example} 
If $\Fq=\seisuu_3[\omega]/(\omega^2+1)$, then $q=3^2$,
$\omega$ has order $8$, and $\nbigq(q)=\{(1,1,3,3)\}$, so
we have a cocycle
$\Gamma(1,1,3,3)=F(1,4,3)-2^{-1}(1-\omega)\cdot F(1,1,6)$.
If $\omega^{a+3}=1$ and $0<a<9$,
then we have $a=1,5$.
Hence we have the cocycles
$$\bigl\{ F(1,3,0),~\Psi(5,3),~\Gamma(1,1,3,3),~E_1(1,3),~E_1(3,9) \bigr\}$$
which form a basis for the third quandle cohomology group.
\end{example}

\begin{example} 
If $\Fq=\seisuu_3[\omega]/(\omega^2+\omega-1)$, then $q=3^2$, $\omega$
has order $8$, and so
$H^3\bigl(\Fq[T]/(T-\omega)\bigr)$ is generated by
the cocycle $\Psi(5,3)$.
\end{example}

\begin{example} 
If $\Fq=\seisuu_2[\omega]/(\omega^3+\omega^2+1)$, then
$q=2^3$ and $\omega$ has order $7$.
We have the triple $(1,2,4)$ of powers of $2$ satisfying $1<2<4$, and 
thus have a cocycle $F(1,2,4)$.
Note that $\omega^{q_1+q_2}\neq 1$ in the case where $q_i$ are powers
of $2$ satisfying $q_i<8$.
Hence the cocycles
$$\bigl\{ F(1,2,4),~ \Psi(5,2),~ \Psi(3,4)\bigr\}$$
form a basis for $H^3\bigl(\Fq[T]/(T-\omega)\bigr)$.
\end{example}

\section{Proof of Theorem~\ref{thm;9.16.1}} \label{section;9.25.3}

\subsection{Preliminaries}
\subsubsection{A decomposition of the complex $C^n$}
We may decompose $C^n$ by the total degree, as follows:
\begin{align*}
 C_d^n&:=\Bigl\{
 \sum a_{i_1,\ldots,i_n}\cdot\prod_{h=1}^{n-1}U_h^{i_h}\cdot T_n^{i_n}
 \in C^n
 ~\Big|~ \sum i_h=d
 \Bigr\}\\
 C_d^n(q)&:=C_d^n\cap C^n(q)
\end{align*}
Then $C^n(q)=\bigoplus_d C_d^n(q)$, and 
it is easy to see that $\delta(C^n_d(q))\subset C_d^{n+1}(q)$.
We denote the complex $\bigl(\bigoplus C^n_d(q),\delta\bigr)$ by
$C^{\ast}_d(q)$.

The following easy lemma follows by a standard argument
(see \cite{m}).
\begin{lem}
In the case $\omega^d\neq 1$, the complex $C^{\ast}_d(q)$ is acyclic.
\end{lem}

The next lemma shows the relationship between this decomposition and
the differential $\delta$.

\begin{lem}
Let $f=\sum_a f_a(U_1,\ldots,U_{n-1})\cdot T_n^a$ be an element of $C^n_d(q)$.
Then
\begin{multline*}
 \delta(f)\bigl(U_1,\ldots,U_{n},T_{n+1}\bigr)=
 \sum_{a}\delta(f_a)\bigl(U_{1},\ldots,U_{n}\bigr)\cdot T_{n+1}^a\\
+(-1)^{n-1}
 \sum_{a}f_a(U_{1},\ldots,U_{n-1})\cdot
 \Bigl(
 \omega^d\cdot\bigl(U_{n}+\omega^{-1}T_{n+1}\bigr)^a
-\bigl(U_{n}+T_{n+1}\bigr)^a
 \Bigr).
\end{multline*}
\end{lem}


\begin{example}
Let $\lambda_d(T_1):=T_1^d\in C^1_d$.  Then
\begin{align*}
 \delta(\lambda_d)(U_1,T_2)&=
(\omega\cdot U_1+T_2)^d-(U_1+T_2)^d\\
&=\omega^d\cdot (U_1+\omega^{-1}\cdot T_2)^d-(U_1+T_2)^d.
\end{align*}
\end{example}
\begin{example}
For an element $f(U_1,T_2)=\sum f_a(U_{1})\cdot T_2^a\in C_d^2$,
\begin{align*}
\delta(f)(U_1,U_2,T_3) &=
 \sum_a \Bigl(
 f_a(\omega U_{1}+U_{2})-f_a(U_{1}+U_2)\Bigr)\cdot T_3^a \\
&\quad\quad-\sum_a f_a(U_{1})\cdot
 \Bigl(
 \omega^d(U_{2}+\omega^{-1}T_3)^a-(U_{2}+T_3)^a\Bigr)\\
&=
 \sum_{a}\delta(f_a)(U_{1},U_{2})\cdot T_3^a\\
&\quad\quad-\sum_af_a(U_{1})\cdot
 \Bigl(
 \omega^d(U_{2}+\omega^{-1}T_3)^a
-(U_{2}+T_3)^a
 \Bigr).
\end{align*}
\end{example}

\begin{example}
If $f=\sum f_a(U_1,U_2)\cdot T_3^a$, then
\begin{align*}
 \delta(f)(U_1,U_2,U_3,T_4)&=
 \sum \delta(f_a)(U_1,U_2,U_3)\cdot T_3^a \\
&\quad\quad+\sum f_a(U_1,U_2)\cdot
 \Bigl(
 \omega^d(U_3+\omega^{-1}T_4)^a
-(U_3+T_4)^a
 \Bigr).
\end{align*}
\end{example}

\subsubsection{The filtration and the derivatives}

Let
\begin{align*}
 C^{n\,(s)}_d&:=
 \Big\{\sum_{a}f_a(U_{1},\ldots,U_{n-1})\cdot T_n^{ap^s}
 \in C^n_d\Big\},\\
C^{n\,(s)}_d(q)&:=C^{n\,(s)}_d\cap C^n_d(q),
\end{align*}
and let
\begin{align*}
  C^{n\,(\infty)}_d&:=
 \Big\{ f_0(U_{1},\ldots,U_{n-1})\in C^n_d(q)\Big\}
=\bigcap_{s\geq 0} C^{n\,(s)}_d(q), \\
C^{n\,(\infty)}_d(q)&:=C^{n\,(\infty)}_d\cap C^n_d(q).
\end{align*}
It is easy to see that $\delta\bigl(C^{n\,(s)}_d(q)\bigr)$
is contained in $C^{n+1\,(s)}_d(q)$, and
we define
\begin{alignat*}{2}
 Z_d^n(q)&:=\Ker(\delta)\cap C_d^n(q), & \qquad
 B_d^n(q)&:=\delta \bigl(C_d^{n-1}(q)\bigr),\\
 Z_d^{n\,(s)}(q)&:=\Ker(\delta)\cap C_d^{n\,(s)}(q),& \qquad
 B_d^{n\,(s)}(q)&:=\delta\bigl(C_d^{n-1\,(s)}(q)\bigr).
\end{alignat*}
There is a homomorphism
$D_n^{(s)}\co C^{n\,(s)}_d(q)\to C^{n\,(s)}_{d-p^s}(q)$
defined as follows:
\[
 D_n^{(s)}\Bigl(\sum_a f_a\cdot T_n^{ap^s}\Bigr)=
\sum_a (a\cdot f_a)\cdot  T_n^{(a-1)p^s}
\]
The kernel $\ker D_n^{(s)} = C^{n\,(s+1)}(q)$, and the
relation $\delta\circ D_n^{(s)}=D_{n+1}^{(s)}\circ\delta$
can be checked easily.
Where the meaning is clear, we may 
omit the subscript $n$.

Let $s$ be a positive integer such that $p^s<q$, and define
\[
 \nbigp(s,q):=
 \bigl\{p^t~|~0\leq t<s\bigr\}
 \cup
 \bigl\{b\cdot p^{s}~|~0<b\cdot p^s<q,~
 b\not\equiv -1~(\modulo p) \mbox{ or } b=p-1 \bigr\}.
\]
For any positive integer $d<q$, set $\lambda_d(T_1):=T_1^d\in C^1_d$.
\begin{lem} \label{lem;04.5.28.20}
Let $s$ and $d$ be integers such that $p^s<q$ and $0<d<q$.
If $\delta(\lambda_d)\in C^2_d(q)$ is contained in
the subset $\Image(D^{(s)})\subset C^{2\,(s)}_d(q)$, then $d\in\nbigp(s,q)$.
\end{lem}
\begin{proof}
We note that $\delta(\lambda_d)=(\omega U_1+T_2)^d-(U_1+T_2)^d$, and
consider the integers $d_t$ such that
$d=\sum d_tp^t$ and $0\leq d_t\leq p-1$ for any non-negative integer $t$.
Set $i:=\min\{t\,|\,d_t>0\}$, then
$(\omega U_1+T_2)^{d-p^i}-(U_1+T_2)^{d-p^i}=0$
for $i<s$, and hence $d-p^i=0$.

If $i\geq s$, then $d$ is of the form $b\cdot p^s$ for some positive
integer $b$.
We will suppose that $b=ap-1$ for some $a>1$, and show that this leads to a
contradiction.
Consider the partition $a=\sum_{t\geq 0} a_t \cdot p^t$
such that $0\leq a_t\leq p-1$.
If $\sum a_t=1$, then $ap=p^h$ for some $h>1$, and so
\[
 \left(\begin{array}{c}a\cdot p-1\\p\end{array}\right)\not\equiv 0
 ~(\modulo p)~\mbox{for}~\omega^{p}-1\neq 0.
\]
Thus the coefficient of $(U_1^p\cdot T_2^{ap-p-1})^{p^s}$
in $\delta(\lambda_d)$ is nonzero.

Now consider the case $\sum a_t>1$, and set $j:=\max\{t~|~a_t>0\}$.
Then
\[
 \left(
 \begin{array}{c}
 a\cdot p-1\\ p^{j+1}
 \end{array}
 \right)
\not\equiv 0 ~(\modulo p)~\mbox{for}~\omega^{p^{j+1}}-1\neq 0,
\]
and so the coefficient of
$\bigl(U^{p^{j+1}}_1T_2^{ap-1-p^{j+1}}\bigr)^{p^s}$
is nonzero, hence
$\delta(\lambda_d)$ cannot be contained in $\Image(D^{(s)})$.
\end{proof}

On the other hand,
if $d\in\nbigp(s,q)$, then $\delta(\lambda_d)$ is contained in
$\Image(D_2^{(s)})$.
For example,
\[
 \delta(\lambda_{(p-1)p^s})
 =
 D_2^{(s)}
 \Bigl[
 p^{-1}\cdot\Bigl(
 (\omega U_1+T_2)^{p}-(U_1+T_2)^p+(1-\omega^p)\cdot U_1^p
 \Bigr)
 \Bigr]^{p^{s}}
\]
when $b=(p-1)p^s$.
Note that $\lambda_{(p-1)p^s}\in \Image(D_2^{(s)})$
even in the case $p^{s+1}=q$.
The other cases can be checked more easily.

\subsubsection{The $2$--cocycles}

A routine calculation proves the following lemma.
\begin{lem}
If $s$ and $t$ are non-negative integers such that
$\omega^{p^s+p^t}=1$, then
$U_1^{p^t}\cdot T_2^{p^s}$ is a quandle 2--cocycle.
\end{lem}

Let $d$ be a positive integer such that $\omega^d=1$, then for any $s$ such
that $p^s<q$, we
consider the following sets of $2$--cocycles:
\[
 J_d^{(s)}(q):=
 \big\{ U_1^{p^t}\cdot T_2^{p^{s}}
 \,\big|\, t<s,\,\,p^s+p^t=d,\,\, p^s<q
 \big\}
\]
\begin{rem}
Clearly the order of $J_d^{(s)}(q)$ is at most $1$.
\end{rem}
Let $H_d^{2\,(s)}(q)$ denote the subspace of $C_d^{2\,(s)}(q)$
generated by $J_d^{(s)}(q)$.

\begin{lem}
If $\omega^d=1$, then
$$Z_d^{2\,(s)}(q)=
 H_d^{2\,(s)}(q)\oplus 
 \bigl(B_d^2(q)\cap C_d^{2\,(s)}(q)+Z_d^{2\,(s+1)}(q)\bigr)$$
and
$Z_d^{2\,(\infty)}(q)=0$.
\end{lem}
\begin{proof}
By an argument similar to
the proof of Lemma \ref{lem;9.11.1} below,
$$Z_d^{2\,(s)}(q)=H_d^{2\,(s)}(q)\oplus
\bigl(B_d^2(q)\cap C_d^{2\,(s)}(q)+Z_d^{2\,(s+1)}(q)\bigr).$$
If $\delta(U_1^d)=0$ for 
$U_1^d\in C^{2\,(\infty)}_d(q)$, then
$d=p^s$ for some $s$, and so
$\omega^d\neq 1$.
\end{proof}


Note that $B_d^2(q)\cap C_d^{2\,(s)}(q)\cap Z_d^{2\,(s+1)}(q)$
is contained in $B^2_d(q)\cap C_d^{2\,(s+1)}(q)$.
Then we obtain, as one of the simplest special cases, the following
proposition, originally stated and proved in \cite{m}.
\begin{prop}
We have a decomposition
$Z_d^2(q)=\bigoplus_{s\geq 0}H_d^{2\,(s)}(q)\oplus B_d^2(q)$.
In particular, the natural map
$\bigoplus H_d^{2\,(s)}(q)\to H^2\bigl(C^{\ast}_d(q)\bigr)$
is an isomorphism.
\end{prop}

\subsubsection{Preliminaries for $3$--coboundaries}

\begin{lem} \label{lem;9.11.1}
If $\omega^d=1$,
then $B_d^{3\,(s)}(q)=B_d^3(q)\cap C_d^{3\,(s)}(q)$
and $B_d^{3\,(\infty)}(q)=B_d^3(q)\cap C_d^{3\,(\infty)}(q)$.
\end{lem}
\begin{proof}
Let $f$ be an element of $C_d^{2\,(s)}$
such that
$f=\sum_a f_{ap^s}(U_{1})\cdot T_2^{ap^s}$.
Then
\[
 \delta(f)=
 \sum_{a} \delta(f_{ap^s})(U_1,U_2)\cdot T_3^{ap^s}
-\sum_a f_{ap^s}(U_1)\cdot
 \Bigl( (U_2+\omega^{-1}T_3)^{ap^s}-(U_2+T_3)^{ap^s}
 \Bigr).
\]
Assume that $\delta(f)\in C_d^{3\,(s+1)}(q)$.
By comparing coefficients of $T_3^{p^s}$,
we find that
\[
 \delta(f_{p^s})(U_1,U_2)
-\sum_a f_{ap^s}(U_1)\cdot a\cdot
 U_2^{(a-1)p^s}\cdot (\omega^{-p^s}-1)=0,
\]
so
$\delta(f_{p^s})+(1-\omega^{-p^s})\cdot D^{(s)}(f)=0$, and hence
$(d-p^s)\in\nbigp(s,q)$.

If $d-p^s=p^t$ for some $t\geq 0$,
then $f(U_1,T_2)=A\cdot U_1^{p^t}\cdot T_2^{p^s}+h$ for some $A\in k$
and $h\in C_d^{2\,(s+1)}(q)$.
If $d-p^s=b\cdot p^s$ for some $b\not\equiv -1\mod p$,
then
$f=A\cdot \delta(\lambda_{d})+h$ for some $A\in k$ and
some $h\in C_d^{2\,(s+1)}(q)$.
We can exclude the possibility $b=p-1$, since in that case
$d=p^s+(p-1)\cdot p^s=p^{s+1}$, and so $\omega^d\neq 1$.

Hence we see that $\delta(f)=\delta(h)$,
and so $B_d^{3\,(s)}(q)\cap C_d^{3\,(s+1)}(q)=B_d^{3\,(s+1)}(q)$, which
implies that
$B_d^{3\,(s)}(q)=B_d^{3}(q)\cap C_d^{3\,(s)}(q)$
for any finite $s$, and also that
$B_d^{3\,(\infty)}(q)=B_d^{3}(q)\cap C_d^{3\,(\infty)}(q)$.
\end{proof}

\subsection{Reductions}

\subsubsection{Subsets of cocycles}

Let $s$ and $t$ be non-negative integers
such that $p^t<p^s<q$.
Then the subsets $I_d^{(s,t)}(q)$
and $I_d^{(s,s)}(q)$
of $I_d(q)$
are defined as follows:
\begin{align*}
 I_d^{(s,t)}(q):=~&
 \bigl\{F(q_1,p^t,p^s)~|~q_1+p^t+p^s=d,\,q_1<p^t\bigr\}
 \\
 &\cup
 \bigl\{\Psi(b\cdot p^t,p^s)~|~
 b\cdot p^t+p^s=d,~b\not\equiv 0~(\modulo p),~b\neq 1
 \bigr\}\\
 &
\cup
 \bigl\{
 E_0(p\cdot p^t,p^s)~|~p^{t+1}+p^s=d
 \bigr\}
 \\
 &\cup
 \bigl\{
 \Gamma(q_1,p^t,q_3,p^s)~|~
 (q_1,p^t,q_3,p^s)\in\nbigq_d(q)
 \bigr\}\\
 I_d^{(s,s)}(q):=~&
 \bigl\{
 E_1(q_1,p\cdot p^s)~|~
 q_1\leq p^s,~q_1+p^{s+1}=d
 \bigr\}
\end{align*}
Set $I_d^{(s)}(q):=\bigcup_{t\leq s}I_d^{(s,t)}(q)$, and
$$I_d^{(\infty)}(q):=
 \bigl\{F(q_1,q_2,0)\,\big|\,
 q_1+q_2=d,\,0<q_1<q_2<q
 \bigr\}.$$
Let $H^{(s,t)}_d(q)$ denote the subspace of $C^{3\,(s)}_d(q)$
generated by $I_d^{(s,t)}(q)$, let
$H^{3\,(s)}_d(q)$ denote $\bigoplus_{t\leq s} H^{(s,t)}_d(q)$, and
let $H^{3\,(\infty)}_d(q)$ denote the subspace
of $C^{3\,(\infty)}_d(q)$ generated by $I_d^{(\infty)}(q)$.

It is easy to see that $I(q)$ and $H^3(q)$ decompose as follows:
\begin{align*}
 I(q)&=
\coprod_d
\Bigl(
I_d^{(\infty)}\sqcup
\coprod_s I_d^{(s)}
\Bigr)\\
H^3(q)&=
\bigoplus_d
\Bigl(
H_d^{3\,(\infty)}(q)
\oplus
\bigoplus_s H^{3\,(s)}_d(q)
\Bigr)
\end{align*}

\subsubsection{First reduction}

The following theorem implies Theorem~\ref{thm;9.16.1}.
\begin{thm} \label{thm;04.5.28.15}
We have the following decomposition:
\[
 Z_d^{3}(q)=
 \Bigl(
 \bigoplus_{s}H_d^{3\,(s)}(q)\oplus H_d^{3\,(\infty)}(q)\Bigr)
 \oplus B_d^{3}(q)
\]
In particular, the natural homomorphism
$$\Bigl(\bigoplus_{s}H_d^{3\,(s)}(q)\Bigr)\oplus H_d^{3\,(\infty)}(q)\to
H^3(C^{\ast}(q))$$
is an isomorphism.
\end{thm}
\begin{proof}
This follows directly from
Lemma~\ref{lem;10.20.5} and Lemma~\ref{lem;9.11.2} below.
\end{proof}

\subsubsection{Second reduction}

\begin{lem} \label{lem;10.20.5}
If $\omega^d=1$, then
$$Z_d^{3\,(\infty)}(q)=
  H_d^{3\,(\infty)}(q)\oplus B_d^{3\,(\infty)}(q).$$
\end{lem}
\begin{proof}
Let $f(U_1,U_2,T_3)=f_0(U_1,U_2)$ be an element of $Z_d^{3\,(\infty)}(q)$.
Then the polynomial $f_0$ is an element of $Z_d^2(q)$.
We have a decomposition $f_0=g+\delta(\lambda_d)$,
where $g$ is an element of $\bigoplus_{s}H^{2\,(s)}_d(q)$,
which gives the required decomposition of $Z_d^{3\,(\infty)}(q)$.
\end{proof}

\begin{lem} \label{lem;9.11.2}
There is a decomposition
$$Z_d^{3\,(s)}(q)=H_d^{3\,(s)}(q)
   \oplus \bigl(B_d^{3\,(s)}(q) + Z_d^{3\,(s+1)}(q)\bigr).$$
\end{lem}
\begin{proof}
This follows from Lemma~3.16 and Lemma~3.17 below.
\end{proof}

\subsubsection{Third reduction}

Let $s$ denote a non-negative integer such that $p^s<q$, and
let $f \in Z_d^{3\,(s)}$ such that
$f = \sum_a f_{ap^s}(U_1,U_2)T_3^{ap^s}$.
Then
\begin{align*}
 \sum_a &\delta(f_{ap^s})(U_1,U_2,U_3) T_4^{ap^s} \\
&+\sum_a f_{ap^s}(U_1,U_2)
 \Bigl(
 (U_3+\omega^{-1}T_4)^{ap^s}
-(U_3+T_4)^{ap^s}
 \Bigr) =0.
\end{align*}
By considering the coefficients of $T_4^{p^s}$,
we find that
\begin{equation} \label{eq;04.5.28.50}
 \delta(f_{p^s})(U_1,U_2,U_3)
+(\omega^{-p^{s}}-1)\cdot
  \sum_a f_{ap^s}(U_1,U_2)\cdot a\cdot U_3^{(a-1)p^s}=0.
\end{equation}
Hence
$\delta(f_{p^s})+(\omega^{-p^{s}}-1)\cdot D_3^{(s)}f=0$, and so
$f_{p^s}$ is contained in $\gyakuzou$,
where 
we denote $D_n^{(s)}\bigl(C_d^{n\,(s)}(q)\bigr)$
by $\Image D_n^{(s)}$ for simplicity.
We thus obtain a map
$\phi\co Z_d^{3\,(s)}(q)\to \gyakuzou$ given by $\phi(f)=f_{p^s}$.
It is clear that $\phi\bigl(Z_d^{3\,(s+1)}(q)\bigr)=0$.

\begin{lem}
If $g\in C_d^{2\,(s)}(q)$ such that $g=\sum g_{ap^s}(U_1)\cdot
T_2^{ap^s}$, then
$$\phi\bigl(\delta(g)\bigr) =
 \delta(g_{p^s})+(1-\omega^{-p^{s}})\cdot D_2^{(s)}(g).$$
\end{lem}
\begin{proof}
This follows by direct calculation.
\end{proof}
We thus obtain a homomorphism:
\[
 \overline{\phi}\co
 \frac{Z_d^{3\,(s)}(q)}{B_d^{3\,(s)}(q)+Z_d^{3\,(s+1)}(q)}
\lrarr
 \frac{\gyakuzou}{\Image D_{2}^{(s)}+B_{d-p^s}^{2\,(s)}(q)}
\]

\begin{lem} \label{lem;10.20.10}
The homomorphism $\overline{\phi}$ is injective.
\end{lem}
\begin{proof}
Let $f \in Z_d^{3\,(s)}(q)$.
If $\overline{\phi}(f)=0$, then $f_{p^s}=\delta(h)+D^{(s)}(g)$
for some $\delta(h)\in B_{d-p^{s}}^{2\,(s)}(q)$ and $g\in C_d^{2\,(s)}(q)$.
Let $\overline{f}$ denote $f-(1-\omega^{-p^s})^{-1}\cdot\delta(g)$.
Then $D^{(s)}(\overline{f})=0$, and so $\overline{f}\in Z_d^{3\,(s+1)}(q)$.
\end{proof}

We denote by $\psi_s$ the composition:
\[
 H_d^{3\,(s)}(q)
 \lrarr
 \frac{Z_d^{3\,(s)}(q)}{B_d^{3\,(s)}(q)+Z_d^{3\,(s+1)}(q)}
\stackrel{\bar{\phi}}{\lrarr}
 \frac{\gyakuzou}{\Image D_2^{(s)}+B_{d-p^s}^{2\,(s)}(q)}
\]

\begin{lem} \label{lem;9.11.3}
The map $\psi_s$ is an isomorphism.
\end{lem}
\begin{proof}
This follows immediately from Proposition~\ref{prop;9.13.1} below.
\end{proof}

\subsubsection{Fourth reduction}
We define
\begin{eqnarray*}
K^{(t)} & := & \delta^{-1}\bigl(\Image(D_3^{(s)})\bigr)\cap
C_{d-p^s}^{2\,(t)}(q) \quad\mbox{for}~ t\leq s,\\
K^{(s+1)} & := & \Image D_2^{(s)}+B^{2\,(s)}_{d-p^s}(q).
\end{eqnarray*}
Then
$$K^{(s+1)}\subset K^{(s)}\subset \cdots \subset K^{(0)}.$$
It can be easily checked that 
the image $\psi_s(H^{3\,(s,t)}_d)$ is contained in $K^{(t)}$
and that $\psi_s(H^{3\,(s,s)}_d)$ is contained in $K^{(s)}$.
There is an induced homomorphism
$$\psi_{(s,t)}\co H_d^{(s,t)}\to K^{(t)}/K^{(t+1)}$$
for all $0\leq t\leq s$.

\begin{prop} \label{prop;9.13.1}
The homomorphisms $\psi_{(s,t)}$ are isomorphisms.
\end{prop}

We now define
{\small
\[
 \nbiga(s,t):=
 \left\{
\begin{array}{ll}
 \left\{(q_1,q_2,q_3)~\left|~
 \begin{array}{l}
 q_1<q_3,~q_2<p^t\leq q_3<p^s,~
  \omega^{q_1+q_3}=1,\\ \omega^{q_1+q_2}\neq 1,
 ~\mbox{if}~\omega^{q_1}=\omega^{q_2}~\mbox{then}~q_1<q_2
 \end{array}\right.
 \right\}
 & \mbox{if}~p\neq 2,\\
\mbox{{}}\\
 \left\{(q_1,q_2,q_3)~\left|~
 \begin{array}{l}
 q_1<q_3,~q_2< p^t\leq q_3\leq p^s,~
  \omega^{q_1+q_3}=1,\\ \omega^{q_1+q_2}\neq 1,
 ~\mbox{if}~\omega^{q_1}=\omega^{q_2}~\mbox{then}~q_1\leq q_2
 \end{array}\right.
 \right\}
 & \mbox{if}~p=2\\
\end{array}
\right.
\] }
and consider the condition
\begin{equation} \label{eq;9.13.2}
 \delta(g)-\sum_{(q_1,q_2,q_3)\in\nbiga(s,t)}
 a_{q_1,q_2,q_3}\cdot U_1^{q_1}\cdot U_2^{q_2}\cdot T_3^{q_3}
 \in \Image(D^{(s)}_3)
\end{equation}
for any element $g\in C^{2\,(t)}_{d-p^s}(q)$.

\begin{prop}\label{prop;9.13.3}
 If there exists an element $g\in C^{2\,(t)}_{d-p^s}(q)$
 satisfying $(\ref{eq;9.13.2})$,
 then all of the coefficients $a_{q_1,q_2,q_3}$ are zero.
\end{prop}

We will prove propositions \ref{prop;9.13.1} and \ref{prop;9.13.3} later,
by descending induction on $t$.
Before going into the proof, we give some remarks:
\begin{itemize}
\item If $p>2$ and $s=t$, then
Proposition $\ref{prop;9.13.3}$ is trivial.
\item If $p>2$, then
either $a_{q_1,q_2,q_3}=0$ or $a_{q_2,q_1,q_3}=0$.
\item If $p=2$ and $q_1\neq q_2$, then
either $a_{q_1,q_2,q_3}=0$ or $a_{q_2,q_1,q_3}=0$.
\end{itemize}

\subsubsection{The case $t=s$}

Let us consider Proposition \ref{prop;9.13.3}
in the case $p=2$ and $s=t$.
By taking the coefficients of $T_3^{p^s}$ in
\eqref{eq;9.13.2},
we see that
\begin{equation}\label{eq:n1}
 \delta(g_{p^s})(U_1,U_2)+(1-\omega^{-p^s})\cdot D_2^{(s)}(g)(U_1,U_2)-
 \sum a_{q_1,q_2,p^s}\cdot
  U_1^{q_1}\cdot U_2^{q_2}=0.
\end{equation}
We have $a_{q_1,q_2,p^s}=0$ unless $q_i<p^s$.
We also know that
either $a_{q_1,q_2,p^s}$ or $a_{q_2,q_1,p^s}$ vanishes if $q_1\neq q_2$,
and so \eqref{eq:n1} has no solution
if one of $a_{q_1,q_2,p^s}$ is nonzero.

Let us consider Proposition~\ref{prop;9.13.1} for general $p$.
We remark that $h(U_1)\cdot T_2^{ap^s}$ is contained
in the image of $D_2^{(s)}$ if $a\not\equiv -1\mod p$, and 
thus we consider
an element $g \in C_{d-p^s}^{2\,(s)}$
of the form
\[
 g(U_1,T_2)=\sum g_{(ap-1)p^s}(U_1)\cdot T_2^{(ap-1)p^s}, 
\]
to obtain
\begin{multline} \label{eq;9.11.5}
 \delta(g)(U_1,U_2,T_3)=
 \sum \delta(g_{(ap-1)p^s})(U_1,U_2)\cdot T_3^{(ap-1)p^s}\\
-\sum g_{(ap-1)p^s}(U_1)\cdot 
 \Bigl(
  \omega^{-p^s}(U_2+\omega^{-1}T_3)^{(ap-1)p^s}-(U_2+T_3)^{(ap-1)p^s}
 \Bigr).
\end{multline}
If $\delta(g)\in \Image(D_3^{(s)})$, then the coefficients in the
right hand side of \eqref{eq;9.11.5} sum to zero.
Taking the terms in $T_3^{(ap-1)p^s}$
and dividing by $T_3^{(p-1)p^s}$,
we see that
\begin{multline}
 \sum \delta(g_{(ap-1)p^s})(U_1,U_2)\cdot T_3^{(a-1)p^{s+1}}\\
-\sum g_{(ap-1)p^s}(U_1)\cdot
 \Bigl(
  \omega^{-p^{s+1}}(U_2+\omega^{-1}T_3)^{(a-1)p^{s+1}}
 -(U_2+T_3)^{(a-1)p^{s+1}}
 \Bigr)=0.
\end{multline}
Substituting $T_3=0$ gives
\[
 \delta(g_{(p-1)p^s})(U_1,U_2)
+  (1-\omega^{-p^{s+1}})\cdot 
 \sum g_{(ap-1)p^s}(U_1)\cdot U_2^{(a-1)p^{s+1}}=0,
\]
which shows that $g(U_1,T_2)=
 (\omega^{-p^s}-1)^{-1}\cdot
  \delta g_{(p-1)p^s}(U_1,T_2)\cdot T_2^{(p-1)p^s}$,
and also that $\delta g_{(p-1)p^s}\in C^{3\,(s+1)}_{d-p^{s+1}}(q)$.
Thus the degree $\deg(g_{(p-1)p^s})=d-p^{s+1}$
is either $p^h$ (for $0\leq h\leq s$) or $b\cdot p^{s+1}$
(see the first half of the proof of Lemma \ref{lem;04.5.28.20}).
In the case $\deg(g_{p^s(p-1)})=b\cdot p^{s+1}$,
the polynomial $g$ is of the form
$\bigl(\omega^{-p^s}-1\bigr)^{-1}\cdot
 \delta(\lambda_{bp^{s+1}})(U_1,T_2)
 \cdot T_2^{(p-1)p^s}$.
Then
\[
 \delta\bigl(\lambda_{(b+1)p^{s+1}-p^s}\bigr)
-\delta\bigl(\lambda_{bp^{s+1}}\bigr)\cdot T_2^{(p-1)p^s}
\in \Image\bigl(D_2^{(s)}\bigr),
\]
and so the term
$\delta\bigl(\lambda_{bp^{s+1}}\bigr)\cdot T_2^{(p-1)p^s}$
can be killed.

On the other hand,
\[
 \psi_s\bigl(E_1(p^{h},p\cdot p^s) \bigr)
=(1-\omega^{-p^s})\cdot
 U_1^{p^h}\cdot T_2^{(p-1)\cdot p^s},
\]
and so the term in $\delta\bigl(\lambda_{p^h}\bigr)\cdot T_2^{(p-1)p^s}$
can be killed by $E_1(p^{h}, p\cdot p^s)$.
Thus we conclude that $\psi_{(s,s)}$ is surjective.


We remark that $p^h+p^s(p-1)\not\equiv 0~(\modulo p^s)$ if $h<s$, and
thus the injectivity of $\psi_{(s,s)}$ can be checked easily.


\subsubsection{The case $t<s$}

We assume that the claims of the propositions \ref{prop;9.13.1}
and \ref{prop;9.13.3} hold for larger than $t+1$,
and we will prove the claims for $t$.

Let $g\in C^{2\,(t)}_{d-p^s}$ be an element
satisfying \eqref{eq;9.13.2}.
Then, comparing
coefficients of $T_3^{p^t}$, we find that
\begin{equation} \label{eq;04.5.28.25}
 \delta(D^{(t)}g)(U_1,U_2,T_3)_{|T_3=0}
=\sum a_{q_1,q_2,p^t}\cdot U_1^{q_1}\cdot U_2^{q_2}.
\end{equation}
Here ``$|T_3=0$'' means the substitution $T_3=0$.
If we decompose $g$ as
$$g=\sum g_{ap^t}(U_1)\cdot T_2^{ap^t},$$
then it is easy to see from \eqref{eq;04.5.28.25} that
$$\delta(g_{p^t})(U_1,U_2)
-\sum a_{q_1,q_2,p^t}\cdot U_1^{q_1}\cdot U_2^{q_2}$$
is contained in $\Image(D_2^{(t)})$.
Recall that we have $a_{q_1,q_2,p^t}=0$ unless $q_i<p^t$.
If $p>2$, or if $p=2$ and $q_1 \neq q_2$, then either
$a_{q_1,q_2,p^t}$ or $a_{q_2,q_1,p^t}$ is zero, and so
we can conclude that all of the $a_{q_1,q_2,p^t}$ are zero
in both cases,
by using
\[
 \delta\bigl(\lambda_{q_1+q_2}\bigr)(U_1,U_2)
=\bigl(\omega^{q_1+q_2}-1\bigr)\cdot
 U_1^{q_1+q_2}
+\bigl(\omega^{q_1}-1\bigr)\cdot U_1^{q_1}\cdot U_2^{q_2}
+\bigl(\omega^{q_2}-1\bigr)\cdot U_1^{q_2}\cdot U_2^{q_1}.
\]
Hence
$\delta(D^{(t)}g)=0$, and so
$D^{(t)}g$ is of the form
\begin{equation} \label{eq;04.5.28.60}
 D^{(t)}g(U_1,T_2)=
 A\cdot \delta(\lambda_{d-p^s-p^t})(U_1,T_2)
+B\cdot U_1^{p_1}\cdot T_2^{p_2}.
\end{equation}
Here we have $p_1<p_2$  and $p^t\leq p_2$.
If $B\neq 0$, then $\omega^{p_1+p_2}=1$.

\begin{lem}
If $p=2$ and $p^t=p_2$, then $B=0$.
\end{lem}
\begin{proof}
Assume that $B\neq 0$.
Since the terms in $D^{(t)}g$ are of the form
$a\cdot U_1^b\cdot T_2^{c\cdot 2^{t+1}}$,
it follows that $A\neq 0$.
We have
\begin{align*}
 \delta\bigl(\lambda_{p_1+p_2}\bigr)
&=\bigl(\omega\cdot U_1+T_2\bigr)^{p_1+p_2}
-\bigl(U_1+T_2\bigr)^{p_1+p_2} \\
&=\bigl(\omega^{p_1}-1\bigr)\cdot U_1^{p_1}\cdot T_2^{p_2}
+\bigl(\omega^{p_2}-1\bigr)\cdot U_1^{p_2}\cdot T_2^{p_1}.
\end{align*}
Note that $p_1<p_2=p^t$, and so
the right hand side of \eqref{eq;04.5.28.60} cannot be contained in 
$\Image(D^{(t)})$, a contradiction.
\end{proof}

Thus $B\cdot U_1^{p_1}\cdot T_2^{p_2}\in \Image(D^{(t)})$.
It follows from Lemma~\ref{lem;04.5.28.20} that
$d-p^s-p^t \in \nbigp(t,q)$
if $A\cdot \delta(\lambda_{d-p^s-p^t})\neq 0$.

\begin{lem}
We can kill
$A\delta(\lambda_{d-p^s-p^t})$
by using one of $D^{(t)}\psi_s(F(p^h,p^t,p^s))$,
$D^{(t)}\psi_s(E_0(p\cdot p^{t},p^s))$ or
$D^{(t)}\psi_s\bigl(\Psi\bigl(b\cdot p^t,p^s\bigr)\bigr)$.
\end{lem}
\begin{proof}
The following may be easily verified by routine calculation.
\begin{align*}
\psi_s\bigl(F(p^h,p^t,p^s)\bigr) &= U_1^{p^h}\cdot T_2^{p^t},\\
\psi_s\bigl(E_0(p\cdot p^t,p^s)\bigr)
 &= \left(\frac{1}{p}\Bigl(\bigl(\omega U_1+T_2\bigr)^p
  -\bigl(U_1+T_2\bigr)^p
  -\bigl(\omega^p-1\bigr)\cdot U_1^p\Bigr)\right)^{p^t},\\
\psi_s\bigl(\Psi(b\cdot p^t,p^s)\bigr)
 &= \bigl(\omega U_1+T_2\bigr)^{b\cdot p^t}
  -\bigl(U_1+T_2\bigr)^{b\cdot p^t}
  -\bigl(\omega^{b\cdot p^t}-1\bigr)\cdot U_1^{b\cdot p^t}.
\end{align*}
As an immediate consequence,
we obtain the following.
\begin{equation} \label{eq;04.5.29.1}
\begin{aligned}
D^{(t)}\psi_s\bigl(F(p^h,p^t,p^s)\bigr)
  &= U_1^{p^h}=(\omega^{p^h}-1)^{-1}\cdot\delta\bigl(\lambda_{p^h}\bigr),\\
 D^{(t)}\psi_s\bigl(E_0(p\cdot p^t,p^s)\bigr)
  &= \bigl(\omega\cdot U_1+T_2\bigr)^{(p-1)\cdot p^t}
  -\bigl(U_1+T_2\bigr)^{(p-1)\cdot p^t} \\
  &=\delta\bigl(\lambda_{p^t\cdot(p-1)}\bigr)(U_1,T_2),\\
 D^{(t)}\psi_s\bigl(\Psi(b\cdot p^t,p^s)\bigr)
  &=b \cdot \Bigl(\bigl(\omega\cdot U_1+T_2\bigr)^{(b-1)\cdot p^t}
  -\bigl(U_1+T_2\bigr)^{(b-1)\cdot p^t}\Bigr) \\
  &=b\cdot \delta\bigl(\lambda_{(b-1)p^t}\bigr).
\end{aligned}
\end{equation}
Then the claim of the lemma follows immediately.
\end{proof}

We now consider the case $B\neq 0$.
\begin{lem}
If $(p_1,p^t,p_2,p^s)\in\nbigq_d(q)$, then
we can kill the term $B\cdot U_1^{p_1}T_2^{p_2}$
by using $D^{(t)}\psi_s(\Gamma(p_1,p^t,p_2,p^s))$.
\end{lem}
\begin{proof}
In cases 1--3 of Condition~\ref{condition;04.5.28.100},
$\psi_s\bigl(\Gamma(p_1,p^t,p_2,p^s)\bigr)
=U_1^{p_1}\cdot T_2^{p^t+p_2}$,
and so
$D^{(t)}\psi_s\bigl(\Gamma(p_1,p^t,p_2,p^s)\bigr)
=U_1^{p_1}\cdot T_2^{p_2}$.
In cases 4 and 5 of Condition~\ref{condition;04.5.28.100},
\begin{multline*}
 \psi_s\bigl(\Gamma(p_1,p^t,p_2,p^s)\bigr)
=U_1^{p_1}\cdot T_2^{p^t+p_2} 
+U_1^{p^t}\cdot T_2^{p_1+p_2} \\
-\bigl(\omega^{q_1}-1\bigr)^{-1}
\cdot \bigl(1-\omega^{2q_1}\bigr)
\cdot U_1^{p_1+p^t}\cdot T_2^{p_2},
\end{multline*}
and so
$$D^{(t)}\psi_s\bigl(\Gamma(p_1,p^t,p_2,p^s)\bigr)
=\left\{
 \begin{array}{ll}
 U_1^{p_1}\cdot T_2^{p_2}, & \mbox{if}~(p^t<p_1),\\
 \mbox{{}}\\
 2\cdot U_1^{p_1}\cdot T_2^{p_2}, & \mbox{if}~(p^t=p_1,\,\,p\neq 2),
 \end{array}
\right.$$
and the claim of the lemma follows immediately.
\end{proof}

\begin{lem}
$g-C\cdot \psi_s(f)$ satisfies (\ref{eq;9.13.2})
for any $f\in H^{(s,t)}_d$ and any $C\in k$.
\end{lem}
\begin{proof}
Since $\delta (f)=0$ for $f\in H^{(s,t)}_d$, we have
$\delta \psi_{s}(f)\in \Image D_2^{(s)}$ due to \eqref{eq;04.5.28.50}, and
the lemma follows immediately.
\end{proof}

\begin{lem}\label{lem;04.5.28.105}
If $(p_1,p^t,p_2,p^s)\not\in\nbigq_d(q)$,
then $(p_1,p^t,p_2)\in\nbiga(s,t+1)$.
\end{lem}
\begin{proof}
We remark that $p_1<p_2$, $p^t<p^{t+1}$, $p^t\leq p_2$
and $\omega^{p_1+p_2}=1$.

The possibility $\omega^{p_1+p^t}=1$
is removed by case 1 of Condition \ref{condition;04.5.28.100}.
We remark that $\omega^{p_1+p^t}\neq 1$ and $\omega^{p_1+p_2}=1$
imply that $p^t\neq p_2$, and so $p^{t+1}\leq p_2$.
If $p\neq 2$, then cases 2 and 3 remove the possibility $p^s\leq p_2$.
Case 4 removes the possibility that both
$p^t\leq p_1$ and $\omega^{p^t}=\omega^{p_1}$ simultaneously, and so we
see that $(p_1,p^t,p_2)\in \nbiga(s,t+1)$ if $p\neq 2$.
If $p=2$, then the claim follows by a similar argument.
\end{proof}

We can choose a constant $C\in k$ such that
we can kill the term $B\cdot U_1^{p_1}\cdot T_2^{p_2}$
by using
$D^{(t)}\bigl(C\cdot U_1^{p_1}\cdot T_2^{p_2+p^t}\bigr)$.
Then $g'=g-C\cdot U_1^{p_1}\cdot T_2^{p_2+p^t}$
satisfies $D^{(t)}(g'-g)=0$,
ie
$g'\in C^{2\,(t+1)}$.
Then we obtain
\begin{align*}
 \delta(g')(U_1,U_2,T_3)&=
 \delta(g)(U_1,U_2,T_3)
-C\cdot\delta(U_1^{p_1}\cdot T_2^{p_2+p^t})\\
 &=\sum_{q_3>p^t} a_{q_1,q_2,q_3}\cdot
  U_1^{q_1} U_2^{q_2} T_3^{q_3}
 -C\cdot (1-\omega^{p_1+p^t})\cdot
  U_1^{p_1} U_2^{p^t} T_3^{p_2}
 +h
\end{align*}
for some $h\in \Image(D_3^{(s)})$.
We conclude that all of the $a_{q_1,q_2,q_3}$ and $C$ are zero,
by Proposition~\ref{prop;9.13.3} for $t+1$
and Lemma~\ref{lem;04.5.28.105}, and so Proposition~\ref{prop;9.13.3} is
true for $t$ as well.

Let $g$ be an element of $K^{(t)}$.
Then $\delta(g)\in \Image\bigl(D^{(s)}_3\bigr)$.
By the argument above,
it can be shown that
$g-\psi_s(f) \in K^{(t+1)}$
for some suitable $f\in H^{(s,t)}_d(q)$.
The linear independence of 
$\bigl\{\psi_s(f)\,\big|\,f\in I^{(s,t)}(q)\bigr\}$
in $K^{(t)}/K^{(t+1)}$ follows
by the argument above, and \eqref{eq;04.5.29.1}.
Thus Proposition~\ref{prop;9.13.1} holds for $t$, and so
the proof of propositions \ref{prop;9.13.1} and \ref{prop;9.13.3}
is complete.

\Addressesr


\begin{thebibliography}{99}

\bibitem{q1}
  \textbf{J\,S Carter}, \textbf{D Jelsovsky}, \textbf{S Kamada},
  \textbf{L Langford}, \textbf{M Saito}, 
  {\em Quandle cohomology and State-sum invariants of knotted curves
  and surfaces},
  Trans. Amer. Math. Soc. 355 (2003) 3947--3989
  \MR{1990571}

\bibitem{q2}
  \textbf{J\,S Carter}, \textbf{D Jelsovsky}, \textbf{S Kamada},
  \textbf{M Saito}, 
  {\em Computations of quandle cocycle invariants of knotted curves
  and surfaces},
  Adv. Math. 157 (2001) 36--94
  \MR{1808844}

\bibitem{q3}
  \textbf{J\,S Carter}, \textbf{D Jelsovsky}, \textbf{S Kamada},
  \textbf{M Saito}, 
  {\em Quandle homology groups, their Betti numbers, and virtual knots},
  J. Pure Appl. Algebra 157 (2001) 135--155
  \MR{1812049}

\bibitem{eg}
 \textbf{P Etingof}, \textbf{M Gra\~{n}a},
 {\em On rack cohomology},
 J. Pure Appl. Algebra 177 (2003) 49--59
  \MR{1948837}

\bibitem{k}
  \textbf{A Kawauchi} (editor),
  {\em A survey of knot theory},
  Birkh\"{a}user,  Basel (1996)
  \MR{1417494}

\bibitem{m}
  \textbf{T Mochizuki},
  {\em Some calculations of cohomology groups of finite
  Alexander quandles},
  J. Pure Appl. Algebra 179, (2003) 287--330
  \MR{1960136}

\bibitem{rs}
  \textbf{C Rourke}, \textbf{B Sanderson},
  {\em A new classification of links and some calculations using it},
  \arxiv{math.GT/0006062}

\bibitem{ss1}
  \textbf{S Satoh}, \textbf{A Shima},
  {\em The $2$--twist spun trefoil has the triple number four},
  Trans. Amer. Math. Soc. 356 (2004) 1007--1024
  \MR{1984465}

\bibitem{ss2}
  \textbf{S Satoh}, \textbf{A Shima},
  {\em Triple point numbers of surface-knots and colorings by
  quandles}, preprint

\end{thebibliography}
\end{document}